\documentclass[12pt]{article}
\usepackage{amsfonts,amsmath,amstext,amssymb}

\newtheorem{thm}{Theorem}
\newtheorem{lemma}[thm]{Lemma}
\newtheorem{proposition}[thm]{Proposition}
\newtheorem{corollary}{Corollary}
\newtheorem{remark}{Remark}
\newtheorem{definition}{Definition}

\newenvironment{demo}{\noindent \textit{Proof.} }{\hfill{$\Box$}}

\title{On the symmetries of a Kaehler manifold}

\author{Alma L. Albujer
\footnote{alma.albujer@uco.es; ORCID: 0000-0002-2384-9061}\\
\and 
Jorge Alc\'azar
\footnote{f82algoj@uco.es}
\and
Magdalena Caballero
\footnote{magdalena.caballero@uco.es; ORCID: 0000-0002-3278-9531}\\
\and
{\small Departamento de Matem\'aticas, Campus de Rabanales,}\\ {\small Universidad de C\'ordoba, 14071 C\'ordoba, Spain}}

\date{}

\begin{document}

\maketitle

\begin{abstract}
	The natural symmetries of Riemannian manifolds are described by the symmetries of its Riemann curvature tensor. In that sense, the most symmetric manifolds are the constant sectional curvature ones. Its natural generalizations are locally symmetric manifolds, semisymmetric manifolds, and pseudosymmetric manifolds. The analogous generalizations of constant holomorphic sectional curvature Kaehler manifolds are locally symmetric Kaehler manifolds, semisymmetric Kaehler manifolds, and holomorphically pseudosymmetric Kaehler manifolds. Do they differ in some way from their Riemannian analogues? Yes, we prove they can all be characterized only in terms of holomorphic planes. Furthermore, the concept of holomorphically pseudosymmetric Kaehler manifold is different from the classical notion of pseudosymmetric Riemannian manifold proposed by Deszcz. We study some relations between both definitions of pseudosymmetry and the so called double sectional curvatures in the sense of Deszcz.  We also present a geometric interpretation of the complex Tachibana tensor and a new characterization of constant holomorphic sectional curvature Kaehler manifolds.

\vspace{.3cm}


\end{abstract}

\section{Introduction}\label{intro}

A Riemannian manifold $(M,g)$ is said to be \emph{locally flat} when its Riemann curvature tensor $R$ vanishes. Along this paper we will consider the Riemann curvature tensor given by
\begin{equation}\label{eq:R}
	R(X,Y)Z=\nabla_X\nabla_Y Z - \nabla_Y\nabla_X Z-\nabla_{[X,Y]}Z,
\end{equation} for any vector fields $X,Y,Z\in\mathfrak{X}(M)$, where $\mathfrak{X}(M)$ stands for the Lie algebra of vector fields on $M$. When necessary, as it is usual in Riemannian geometry, we will also consider $R$ as the $(0,4)$-tensor given by $R(X,Y,Z,W)=g(R(X,Y)Z,W)$ for any $X,Y,Z,W\in\mathfrak{X}(M)$.

The simplest non-flat Riemannian manifolds are the spaces of constant sectional curvature $c$, for which the Riemann curvature tensor takes the form 
\begin{equation}\label{eq:CSC}
R(X,Y)Z=c(X\wedge_g Y)Z,
	\end{equation}
$\wedge_g$ being the metric endomorphism in $(M,g)$, defined by $X\wedge_g Y:=g(Y,\cdot)X-g(X,\cdot)Y$, for any vector fields $X,Y \in \mathfrak{X}(M)$. Generalizing those spaces, Cartan considered in~\cite{Cartan} the \emph{locally symmetric} spaces as those such that $\nabla R=0$, that is, with parallel curvature, characterized by the fact that the sectional curvature of every plane is preserved after parallel transport of the plane along any curve, see~\cite{Levy}.

Going a step further, let us observe that the integrability condition of $\nabla R=0$ is $R\cdot R =0$. Therefore, every locally symmetric space also satifies $R\cdot R=0$. This fact led Szab\'o~(\cite{Szabo_local,Szabo_global}) to define \emph{semisymmetric spaces} as the Riemannian manifolds satisfying $R\cdot R=0$. Let us recall here that $R\cdot R$ is the $(0,6)$-tensor on $M$ given by
\begin{equation}\label{eq:RR}
	\begin{split}
	R\cdot R(X_1,X_2,X_3,X_4;X,Y)=(R(X,Y)\cdot R)(X_1,X_2,X_3,X_4)\\
	=-R(R(X,Y)X_1,X_2,X_3,X_4)-R(X_1,R(X,Y)X_2,X_3,X_4)\\
	-R(X_1,X_2,R(X,Y)X_3,X_4)-R(X_1,X_2,X_3,R(X,Y)X_4),
\end{split}
\end{equation}
for every $X_1,X_2,X_3,X_4,X,Y\in\mathfrak{X}(M)$. As proved by Haesen and Verstraelen in~\cite[Corollary 1]{HV_2007}, semisymmetric manifolds are characterized as those manifolds whose sectional curvature at every point $p\in M$ of any plane $\pi\subset T_pM$, $K(p,\pi)$, is invariant, up to second order, under parallel transport of $\pi$ around any infinitesimal coordinate parallelogram centered at $p$. It is immediate to observe that any Riemannian surface satisfies $R\cdot R=0$, so the concept of semisymmetry takes relevance in the case of Riemannian manifolds $(M^n,g)$ of dimension $n\geq 3$.

As a natural generalization of these last spaces, the \emph{pseudosymmetric spaces in the sense of Deszcz} (see~\cite{Deszcz} and references therein) are defined as the Riemannian manifolds $(M^n,g)$ of dimension $n\geq 3$ for which
\begin{equation}\label{eq:ps_D}
	R\cdot R=L \,\,Q(g,R),
\end{equation}
where $L\in\mathcal{C}^\infty(M)$ and $Q(g,R)$ is the \emph{Tachibana tensor}, which is defined by
\begin{equation}\label{eq:Tach_real}
		\begin{split}
		Q(g,R)(X_1,X_2,X_3,X_4;X,Y)=-((X\wedge_g Y)\cdot R)(X_1,X_2,X_3,X_4)\\
		=R((X\wedge_g Y)X_1,X_2,X_3,X_4)+R(X_1,(X\wedge_g Y)X_2,X_3,X_4)\\
		+R(X_1,X_2,(X\wedge_g Y)X_3,X_4)+R(X_1,X_2,X_3,(X\wedge_g Y)X_4),
	\end{split}
\end{equation}
for $X_1,X_2,X_3,X_4,X,Y\in\mathfrak{X}(M)$. The Tachibana tensor may well be the simplest $(0,6)$-tensor with the same algebraic symmetries as $R\cdot R$.

It is well known that a Riemannian manifold has constant sectional curvature if and only if its Tachibana tensor vanishes identically, see~\cite{Eisenhart}. Therefore, if $(M^n,g)$, $n\geq 3$, is a Riemannian manifold with non-constant sectional curvature, the set of points where the Tachibana tensor does not vanish identically is an open non-empty subset $\mathcal{U}\subseteq M$. Given a point $p\in\mathcal{U}$, a plane $\pi=v\wedge w\subset T_pM$ is said to be \emph{curvature-dependent} with respect to another plane $\bar{\pi}=x\wedge y\subset T_pM$ if $Q(g,R)(v,w,w,v;x,y)\neq 0$. In this context, given $p\in \mathcal{U}$ and two planes in $T_pM$, $\pi=v\wedge w$ and $\bar{\pi}=x\wedge y$, such that $\pi$ is curvature-dependent with respect to $\bar{\pi}$, the \emph{double sectional curvature} $L(p,\pi,\bar{\pi})$ of $\pi$ with respect to $\bar{\pi}$ at $p$ is defined in~\cite[Definition 3]{HV_2007} as
\begin{equation}\label{eq:dousec}
L(p,\pi,\bar{\pi})=\frac{R\cdot R(v,w,w,v;x,y)}{Q(g,R)(v,w,w,v;x,y)}.
\end{equation}
The above definitions are independent on the choice of bases for $\pi$ and $\bar{\pi}$.

Given $p\in M$ and two planes, $\pi=v\wedge w$ and $\bar{\pi}=x\wedge y$, such that $\pi$ is curvature-dependet with respect to $\bar{\pi}$, $Q(g,R)(v,w,w,v;x,y)$ measures the change of the sectional curvature $K(p,\pi)$ under an operation involving infinitesimal rotations performed at the point $p$. Specifically, considering $\pi'=v'\wedge w'$, where $v',w'\in T_pM$ are vector obtained after infinitesimal rotations of the projections of $v$ and $w$ respectively onto $\bar{\pi}$, $Q(g,R)(v,w,w,v;x,y)$ measures, up to order two, the difference $K(p,\pi')-K(p,\pi)$, see~\cite{HV_2007}.

It is obvious that if $(M,g)$ is a pseudosymmetric manifold in the sense of Deszcz, all its double sectional curvatures are independent of the planes. Furthermore, this is also a sufficient condition: a Riemannian manifold $M$ of dimension $n\geq 3$ is pseudosymmetric in the sense of Deszcz if and only if at every $p\in\mathcal{U}$, for all planes $\pi$ and $\bar{\pi}$ in $T_p M$ such that $\pi$ is curvature-dependent with respect to $\bar{\pi}$, $L(p,\pi,\bar\pi)=L(p)$ for some smooth function $L\in\mathcal{C}^\infty(M)$ (see~\cite[Theorem 3]{HV_2007}). Moreover, $M$ is pseudosymmetric in the sense of Deszcz if and only if at every $p\in\mathcal{U}$, for every planes $\pi$ and $\bar{\pi}$ in $T_p M$ such that $\pi$ is curvature-dependent with respect to $\bar{\pi}$, the double sectional curvature $L(p,\pi,\bar\pi)$ is independent of the plane $\pi$ (see~\cite[Theorem 5]{HV_2007}).

\begin{table}[h]
	\begin{center}
		\begin{tabular}{ccc}
			\hline
			Locally flat & $\qquad$ & $R=0$\\
			\hline
			Constant sectional curvature & & $R(X,Y)Z=c\,(X\wedge_g Y)Z$\\
			\hline
			Locally symmetric & & $\nabla R =0$\\
			\hline
			Semisymmetric & & $R\cdot R =0$\\
			\hline
			Deszcz pseudosymmetric & & $R\cdot R = L\,Q(g,R)$\\
			\hline
		\end{tabular}
		\caption{Natural symmetries on a Riemannian manifold.}
	\end{center}
\end{table}

Along this work, we will consider Kaehler manifolds $(M^{2n},g,J)$ satisfying different symmetries. Specifically, we will study semisymmetric and pseudosymmetric Kaehler manifolds. It is worth pointing out at this point that every Kaehler manifold $(M^{2n},g,J)$ of real dimension $2n>4$, which is pseudosymmetric in the sense of Deszcz, is also semisymmetric, see~\cite{DDV} and~\cite{Olszak2003}. However, there exists an example in dimension $2n=4$ of a Kaehler manifold $(M^4,g,J)$ which is pseudosymmetric in the sense of Deszcz but not semisymmetric,~\cite{Olszak2003}. Furthermore, the tensors $R\cdot R$ and $Q(g,R)$ do not present the same symmetries and properties involving the complex structure $J$. For these reasons, in 1989 Olszak,~\cite{Olszak89}, proposed an alternative concept of pseudosymmetry for Kaehler manifolds. We will refer to this concept of pseudosymmetry as \emph{holomorphic pseudosymmetry}.
	
The manuscript is organized as follows. In Section~\ref{sec:setup}, some basic notions about Kaehler manifolds 
are recalled and three new algebraic results concerning the symmetries of $(0,5)$ and $(0,6)$-tensors on $M$ are presented: Lemma~\ref{lem:simRRJ}, Proposition~\ref{prop:auxalg}, and Proposition~\ref{prop:auxalg2}. In Section~\ref{sec:tachibana}, the complex Tachibana tensor is defined and, in an analogous way as in the Riemannian case, a geometrical interpretation is given. Such tensor allows us to give a characterization result of Kaehler manifolds with constant holomorphic sectional curvature in terms of holomorphic planes: Theorem~\ref{thm:chsc}. In Section~\ref{sec:ls_ss}, characterization results for locally symmetric and semisymmetric Kaehler manifolds are proved: Theorems~\ref{thm:locsym-charac} and~\ref{thm:semi-charac}. Finally,  Section~\ref{sec:holom_ps} is devoted to holomorphically pseudosymmetric Kaehler manifolds. The relations between holomorphic pseudosymmetry, pseudosymmetry in the sense of Deszcz and the double sectional curvatures are studied. In particular, two characterization results for holomorphically pseudosymmetric Kaehler manifolds are given in terms of the double sectional curvatures: Theorems~\ref{thm:DvsO} and~\ref{thm:DvsO_weak}.

\section{Set up}\label{sec:setup}

A Kaehler manifold is a complex manifold $(M^{2n},g,J)$, of real dimension $2n$, where $g$ is an Hermitian metric on $M$ and the complex structure $J$ is parallel, i.e. $\nabla J=0$. Given $(M,g,J)$ a Kaehler manifold, it is possible to define its complex metric endomorphism by
\begin{equation}\label{eq:cmh}
	(X\wedge_g^c Y)Z=(X\wedge_g Y)Z+(JX\wedge_g JY)Z-2g(JX,Y)JZ,
\end{equation}
for any $X,Y,Z\in\mathfrak{X}(M)$. The simplest non-flat Kaehler manifolds are those of constant holomorphic sectional curvature $\tilde{c}$, for which the Riemann curvature tensor takes the form
\begin{equation}\label{eq:CHSC}
	R(X,Y)Z=\dfrac{\tilde{c}}{4}(X\wedge_g^c Y)Z,
\end{equation}
for any $X,Y,Z\in\mathfrak{X}(M)$.
	
The following lemma, due to Ogiue~\cite{Ogiue}, gives a sufficient and necessary condition for a Kaehler manifold to have constant holomorphic sectional curvature.
	
\begin{lemma}\label{lem:Ogiue}
	A Kaehler manifold $(M,g,J)$ has constant holomorphic sectional curvature if and only if $R(X,JX,X,Y)=0$ for every orthonormal subset $\{X,JX,Y\}$ in $\mathfrak{X}(M)$.
\end{lemma}

When studying further symmetries of any Kaehler manifold, or in general of any Riemannian manifold, the tensor $R\cdot R$ proves to be essential. In addition to the algebraic symmetries satisfied by such tensor on any Riemannian manifold (see, for instance,~\cite[Lemma 1]{HV_2007}), in the case of a Kaehler manifold it presents a good behaviour with respect to $J$. Specifically, we show the following lemma.

\begin{lemma}\label{lem:simRRJ}
	The tensor $R\cdot R$ of a Kaehler manifold $(M,g,J)$ verifies
	\begin{equation}\label{eq:simRRJ}
		\begin{split}
			&R\cdot R(JX_1,JX_2,X_3,X_4;X,Y)=R\cdot R(X_1,X_2,JX_3,JX_4;X,Y)\\
			&=R\cdot R(X_1,X_2,X_3,X_4;JX,JY)=R\cdot R(X_1,X_2,X_3,X_4;X,Y),
		\end{split}
	\end{equation}
	for every $X_1,X_2,X_3,X_4,X,Y\in\mathfrak{X}(M)$.
\end{lemma}
\begin{demo}
	The proof follows immediately from the definition of $R\cdot R$ and taking into account~\cite[Prop. 4.5, Ch. IX]{KN_1969}, which asserts that the curvature tensor of a Kaehler manifold satisfies
	\begin{equation}\label{eq:simRJ}
		R(JX,JY)Z=R(X,Y)Z \quad \textrm{and} \quad R(X,Y)JZ=J(R(X,Y)Z),
	\end{equation}
	for every $X,Y,Z\in\mathfrak{X}(M)$.
\end{demo}

Thanks to the following result, $(0,6)$-tensors on a Kaehler manifold satisfying the same properties as $R\cdot R$ are characterized by their behaviour when applied to holomorphic planes.

\begin{proposition}\label{prop:auxalg}
	Let $V^{2n}$ be a real vectorial space endowed with a complex structure $J$ 
	and let $T_1$ and $T_2$ be two $(0,6)$-tensors on $V$ satisfying the following properties:
	\begin{itemize}
		\item[$(a)$] $T_i(x_1,x_2,x_3,x_4,x_5,x_6)=-T_i(x_2,x_1,x_3,x_4,x_5,x_6)=-T_i(x_1,x_2,x_4,x_3,x_5,x_6)$
		$\hspace*{3.77cm}=T_i(x_3,x_4,x_1,x_2,x_5,x_6)$,
		\item[$(b)$] $T_i(x_1,x_2,x_3,x_4,x_5,x_6)+T_i(x_1,x_3,x_4,x_2,x_5,x_6)+T_i(x_1,x_4,x_2,x_3,x_5,x_6)=0$,
		\item[$(c)$] $T_i(x_1,x_2,x_3,x_4,x_5,x_6)=T_i(Jx_1,Jx_2,x_3,x_4,x_5,x_6)=T_i(x_1,x_2,Jx_3,Jx_4,x_5,x_6)$, and
		\item[$(d)$] $T_i(x_1,x_2,x_3,x_4,x_5,x_6)=-T_i(x_1,x_2,x_3,x_4,x_6,x_5)=T_i(x_1,x_2,x_3,x_4,Jx_5,Jx_6)$,
	\end{itemize}
	for every $x_1,x_2,x_3,x_4,x_5,x_6 \in V$. If for every $u,v\in V$ it holds
	\begin{equation}\label{eq:propauxalg}
		T_1(u,Ju,Ju,u,v,Jv)=T_2(u,Ju,Ju,u,v,Jv),
	\end{equation}
	then $T_1=T_2$.
\end{proposition}
\begin{demo}
	Observe first that we can assume without loss of generality that $T_2=0$.
	
	Given $v\in V$, let us consider the $(0,4)$-tensor on $V$, $T^v$, defined by 
	\begin{equation}\label{eq:Tv}
		T^v(x_1,x_2,x_3,x_4)=T_1(x_1,x_2,x_3,x_4,v,Jv).
	\end{equation} 
	From the symmetries of $T_1$ we can immediately check that $T^v$ satisfies
	\begin{itemize}
		\item[$(a^v)$] $T^v(x_1,x_2,x_3,x_4)=-T^v(x_2,x_1,x_3,x_4)=-T^v(x_1,x_2,x_4,x_3)=T^v(x_3,x_4,x_1,x_2)$,
		\item[$(b^v)$] $T^v(x_1,x_2,x_3,x_4)+T^v(x_1,x_3,x_4,x_2)+T^v(x_1,x_4,x_2,x_3)=0$, and
		\item[$(c^v)$] $T^v(x_1,x_2,x_3,x_4)=T^v(Jx_1,Jx_2,x_3,x_4)=T^v(x_1,x_2,Jx_3,Jx_4)$,
	\end{itemize}
	for every $x_1,x_2,x_3,x_4\in V$. Furthermore, from~\eqref{eq:propauxalg} it also holds
	\begin{equation}\label{eq:propauxalgv}
		T^v(u,Ju,Ju,u)=0 \quad \textrm{for every } u\in V.
	\end{equation}
	Therefore, from~\cite[Prop. 7.1, Ch. IX]{KN_1969} we get $T^v=0$, i.e.
	\begin{equation}\label{eq:propauxalgT1}
		T_1(x_1,x_2,x_3,x_4,v,Jv)=0 \quad \textrm{for every } x_1,x_2,x_3,x_4,v\in V.
	\end{equation}
	
	Given now $x_5,x_6\in V$, let us consider $v=Jx_5+x_6$. Then, $Jv=-x_5+Jx_6$ and~\eqref{eq:propauxalgT1} yields
	\begin{equation}\label{eq:propauxalgT12}
		0=T_1(x_1,x_2,x_3,x_4,v,Jv)=T_1(x_1,x_2,x_3,x_4,Jx_5,Jx_6)+T_1(x_1,x_2,x_3,x_4,x_6,-x_5).
	\end{equation}
	Finally, the result follows from assumption $(d)$.
	
\end{demo}

\vspace{.5cm}

\noindent Our last algebraic result is proven in an analogous way.

\begin{proposition}\label{prop:auxalg2}
	Let $V^{2n}$ be a real vectorial space endowed with a complex structure $J$ 
	and let $T_1$ and $T_2$ be two $(0,5)$-tensors on $V$ satisfying the following properties:
	\begin{itemize}
		\item[$(a)$] $T_i(x_1,x_2,x_3,x_4,x_5)=-T_i(x_2,x_1,x_3,x_4,x_5)=-T_i(x_1,x_2,x_4,x_3,x_5)=T_i(x_3,x_4,x_1,x_2,x_5)$,
		\item[$(b)$] $T_i(x_1,x_2,x_3,x_4,x_5)+T_i(x_1,x_3,x_4,x_2,x_5)+T_i(x_1,x_4,x_2,x_3,x_5)=0$, and
		\item[$(c)$] $T_i(x_1,x_2,x_3,x_4,x_5)=T_i(Jx_1,Jx_2,x_3,x_4,x_5)=T_i(x_1,x_2,Jx_3,Jx_4,x_5)$,
	\end{itemize}
	for every $x_1,x_2,x_3,x_4,x_5 \in V$. If for every $u,v\in V$ it holds
	\begin{equation}\label{eq:propauxalg2}
	T_1(u,Ju,Ju,u,v)=T_2(u,Ju,Ju,u,v),
	\end{equation}
	then $T_1=T_2$.
\end{proposition}

\section{The complex Tachibana tensor and constant holomorphic sectional curvature Kaehler manifolds}\label{sec:tachibana}

In analogy to the classical Tachibana tensor, on any Kaehler manifold we define the complex Tachibana tensor as follows.

\begin{definition}
	Given $(M^{2n},g,J)$ a Kaehler manifold, the complex Tachibana tensor is defined as the $(0,6)$-tensor on $M$ given by
	\begin{equation}\label{eq:TCdef}
		\begin{split}
		Q^c(g,R)(X_1,X_2,X_3,X_4;X,Y)=-\left(\left(X\wedge_g^c Y\right)\cdot R\right)(X_1,X_2,X_3,X_4)\\
		=R\left(\left(X\wedge_g^c Y\right)X_1,X_2,X_3,X_4\right)+ R\left(X_1,\left(X\wedge_g^c Y\right)X_2,X_3,X_4\right)\\
		+R\left(X_1,X_2,\left(X\wedge_g^c Y\right)X_3,X_4\right)+ R\left(X_1,X_2,X_3,\left(X\wedge_g^c Y\right)X_4\right),
		\end{split}
	\end{equation}
for every $X_1,X_2,X_3,X_4,X,Y\in\mathfrak{X}(M)$.
\end{definition}

\begin{remark}\label{rmk:Pi}
The choice of the notation $Q^c(g,R)$ for the complex Tachibana tensor is due to its analogy with the classical Tachibana tensor for a Riemannian manifold. However, variations of such tensor have been previously considered in the study of the symmetries of a Kaehler manifold. Specifically, Olszak considered in~\cite{Olszak89} 
the $(0,6)$-tensor $R_1\cdot R$ given by $R_1\cdot R=-Q^c(g,R)$ and Jelonek defined in~\cite{Jelonek} the tensor $\Pi\cdot R=-\frac{1}{4} Q^c(g,R)$. Here $\Pi$ stands for the tensor on a Kaehler manifold given by $\Pi(X,Y)=\frac{1}{4}(X\wedge_g^c Y)$ for any $X,Y\in\mathfrak{X}(M)$. Sometimes it will be convenient to consider $\Pi$ instead of the complex metric endomorphism.
\end{remark}

Since the complex metric endomorphism is, up to a constant, the curvature tensor of a complex space form, the complex Tachibana tensor may well be the simplest non-trivial $(0,6)$-tensor with the same algebraic symmetries and the same properties with respect to $J$ as $R\cdot R$. In fact, one could think that the simplest tensor with the same properties as $R\cdot R$ would be $\Pi\cdot \Pi$. However, with a tedious but straightforward computation it is possible to check that $\Pi\cdot \Pi=0$. Furthermore, taking into account~\eqref{eq:cmh}, for every $p\in M$ and every $v,w,x,y\in T_pM$, it is immediate to get the following relation between the complex and the classical Tachibana tensor,
\begin{equation}\label{eq:TC}
\begin{split}
Q^c(g,R)(v,w,w,v;x,y)=&Q(g,R)(v,w,w,v;x,y)+Q(g,R)(v,w,w,v;Jx,Jy)\\&-4 g(Jx,y)R(Jv,w,w,v)-4 g(Jx,y)R(v,Jw,w,v)\\=&Q(g,R)(v,w,w,v;x,y)+Q(g,R)(v,w,w,v;Jx,Jy).\end{split}
\end{equation} 

Given $p\in \mathcal{U}$ and two planes in $T_pM$, $\pi=v\wedge w$ and $\bar{\pi}=x\wedge y$, it is easy to check that $Q(g,R)(v,w,w,v;x,y)$ and $Q^c(g,R)(v,w,w,v;x,y)$ do not depend on the choice of bases for $\pi$ and $\bar{\pi}$.  And so, we can simply write $Q(g,R)(\pi;\bar{\pi})$ and $Q^c(g,R)(\pi;\bar{\pi})$. Taking into account this notation,~\eqref{eq:TC} reads
\begin{equation}\label{eq:TC2}
Q^c(g,R)(\pi;\bar{\pi})=Q(g,R)(\pi;\bar{\pi})+Q(g,R)(\pi;J\bar{\pi}),
\end{equation}
where $J\bar\pi=Jx\wedge Jy$. It is also worth pointing out that 	whenever any of the two planes is a complex holomorphic plane, i.e. invariant under the action of $J$, the Tachibana tensor and the complex Tachibana tensor are proportional. Indeed, it is easy to check that 
\begin{equation}\label{eq:TQ1}
Q^c(g,R)(\pi^h;\bar{\pi})=2 Q(g,R)(\pi^h;\bar{\pi}),
\end{equation}
\begin{equation}\label{eq:TQ2}
Q^c(g,R)(\pi;\bar{\pi}^h)=2 Q(g,R)(\pi;\bar{\pi}^h),
\end{equation}
where we use the superindex $^h$ for holomorphic planes.

As a first result in this section, we get the following characterization for Kaehler manifolds with constant holomorphic sectional curvature.

\begin{thm}\label{thm:chsc}
	Let $(M,g,J)$ be a Kaehler manifold, then the following conditions are equivalent:
	\begin{itemize}
		\item[$(a)$] $M$ has constant holomorphic sectional curvature.
		\item[$(b)$] The complex Tachibana tensor $Q^c(g,R)$ vanishes identically.
		\item[$(c)$] $Q^c(g,R)(u,Ju,Ju,u;x,Jx)=0$ for every $p\in M$ and $x,u\in T_pM$.
		\item[$(d)$] $Q(g,R)(u,Ju,Ju,u;x,Jx)=0$ for every $p\in M$ and $x,u\in T_pM$.
	\end{itemize}
\end{thm}
\begin{demo}
	Implications $(b) \Rightarrow (c) \Rightarrow (d)$ are inmmediate. For the proof of the implication $(a) \Rightarrow (b)$, just observe that if $M$ has constant holomorphic sectional curvature, its Riemann curvature tensor takes the form $R=c\,\Pi$, where $\Pi$ is defined as in Remark~\ref{rmk:Pi}. Then $Q^c(g,R)=-4\Pi\cdot R=-4c\,\Pi\cdot\Pi=0$.
	
	It remains to show $(d) \Rightarrow (a)$, which is a consequence of Lemma~\ref{lem:Ogiue}.  Given $p\in M$ and $x,u\in T_pM$,
	$$0=Q(g,R)(u,Ju,Ju,u;x,Jx)$$ $$=4g(u, x)R(u, x, Ju, u) - 4g(Ju, x)R(u, Jx, Ju, u).$$
	And so, $$g(u, x)R(u, x, Ju, u) = g(Ju, x)R(u, Jx, Ju, u).$$
	
	We decompose $x$ in its projection onto the plane generated by $u$ and $Ju$ and its orthogonal projection, $x=\alpha u+ \beta Ju+w$, and we get
	$$\alpha g(u,u)\{\beta R(u, Ju, Ju, u)+R(u, w, Ju, u)\}$$ $$= \beta g(Ju,Ju)\{\alpha R(u, Ju, Ju, u)+R(u, Jw, Ju, u)\}.$$
	And so, 
	$$\alpha R(u, w, Ju, u) = \beta R(u, Jw, Ju, u).$$
	In particular, if $\beta$ vanishes and $\alpha\neq 0$, $$R(u, w, Ju, u)=0.$$
	
	We choose now an orthogonal set in $T_pM$, $\{u,Ju,w\}$, and we construct $x=u+w$. Applying the reasoning above, we finish the proof taking into account Lemma~\ref{lem:Ogiue}.
\end{demo}		

Our next aim is to give a geometrical interpretation of the tensor $Q^c(g,R)$. As a previous step, we recall the geometrical interpretation of the metric endomorphism applied to a vector. Even if the interpretation is valid for any Riemannian manifold, we focus on Kaehler manifols.   

Let $(M,g,J)$ be a Kaehler manifold. Take $p\in \mathcal{U}$, a plane $\bar{\pi}=x\wedge y$ in $T_pM$ and $z \in T_pM$. Assume that $x,y$ are orthonormal and choose $\{e_3,\cdots, e_{2n}\}$, so that $\{x,y,e_3,\cdots, e_{2n}\}$ is an orthonormal basis of $T_p M$. Then $z$ can be decomposed as the sum of its projection onto $\bar{\pi}$ and its projection onto the $(2n-2)$-linear subspace of $T_pM$ expanned by $e_3,\cdots, e_{2n}$. By rotating the projection of $z$ onto $\bar{\pi}$ an angle $\varepsilon$, while keeping the other part of the sum fixed, a new vector, $\tilde{z}$, is obtained. If $\varepsilon$ is small enough, that vector can be aproximated by 

\begin{equation}\label{eq:IntGeom1}
\tilde{z}=z+\varepsilon (x\wedge_g y)z+O(\varepsilon^2).
\end{equation}

And so, the vector $(x\wedge_g y)z$ measures the first order change of the vector $z$ after such an infinitesimal rotation of $z$ in the plane $x\wedge y$ at the point $p$.  

Take now $p\in \mathcal{U}$ and two planes in $T_pM$, $\pi=v\wedge w$ and $\bar{\pi}=x\wedge y$, where $\{v,w\}$ as well as $\{x,y\}$ are orthonormal. We rotate infinitesimally the projection of $v$ and $w$ onto $\bar{\pi}$. Proceeding as in \eqref{eq:IntGeom1}, we get 
\begin{equation*}\label{eq:IntGeom2}
\begin{split}\tilde{v}=&v+\varepsilon (x\wedge_g y)v+O(\varepsilon^2),\\
\tilde{w}=&w+\varepsilon (x\wedge_g y)w+O(\varepsilon^2).
\end{split}
\end{equation*}

After that, we rotate infinitesimally the projection of those new vectors onto $J\bar{\pi}=Jx\wedge Jy$ to obtain
\begin{equation*}\label{eq:IntGeom3}
\begin{split}
\tilde{v}'=&v+\varepsilon (x\wedge_g y)v+\varepsilon (Jx\wedge_g Jy)v+O(\varepsilon^2),\\
\tilde{w}'=&w+\varepsilon (x\wedge_g y)w+\varepsilon (Jx\wedge_g Jy)w+O(\varepsilon^2).
\end{split}
\end{equation*}

Comparing the sectional curvatures of the planes $\pi=v\wedge w$ and $\tilde{\pi}'=\tilde{v}'\wedge \tilde{w}'$, we get
\begin{equation}\label{eq:IntGeom4}
\begin{split}
K(p, \tilde{\pi}')=&K(p,\pi)+\varepsilon Q(g,R)(\pi,\bar{\pi})+\varepsilon Q(g,R)(\pi,J\bar{\pi})+O(\varepsilon^2)\\
=&K(p,\pi)+\varepsilon Q^c(g,R)(\pi;\bar{\pi})+O(\varepsilon^2).
\end{split}
\end{equation}
In conclusion, the complex Tachibana tensor $Q^c(g,R)(\pi;\bar{\pi})$ measures the change of the sectional curvature of $\pi$ under an operation involving infinitesimal rotations in $\bar{\pi}$ and $J\bar{\pi}$.

\section{Locally symmetric and semisymmetric Kaehler manifolds}\label{sec:ls_ss}

As we recalled in the introduction, locally symmetric spaces, i.e. those for which $\nabla R=0$, are a natural generalization of spaces of constant curvature. 
As we have pointed out in the introduction, it is well known that those spaces are characterized by the fact that the sectional curvature is conserved under parallel transport along any curve. Our next result for Kaehler manifols shows that in the previous characterization we can only consider holomorphic planes.
 
\begin{thm}\label{thm:locsym-charac}
	Let $(M,g,J)$ be a Kaehler manifold. Then, the following conditions are equivalent:
	\begin{itemize}
		\item[$(a)$] $M$ is locally symmetric.
		\item[$(b)$] $\left(\nabla_X R\right)(U,JU,JU,U)=0$ for any vector fields $X,U$ on $M$.
		\item[$(c)$] 
		Sectional curvature of holomorphic planes is invariant under parallel transport along any curve. 
	\end{itemize}
\end{thm}

\begin{demo}
	The proof of $(a)\Rightarrow(b)$ is trivial and  $(a)\Rightarrow(c)$ 
	is known. To finish the proof we will prove $(b)\Rightarrow(a)$ and $(c)\Rightarrow(a)$.

	The implication $(b)\Rightarrow(a)$ is obtained by applying Proposition \ref{prop:auxalg2} to the $(0,5)$-tensor given by $T(X_1,X_2,X_3,X_4,X_5):=\left(\nabla_{X_5}R\right)(X_1,X_2,X_3,X_4)$.
	
	In~\cite[Prop. 10, Ch. 8]{ONeill}, it is proven that being locally symmetric is equivalent to the following property: if $X, Y, Z$ are parallel vector fields along a curve $\alpha$ on $M$, so it is $R(X,Y,Z)$. We use this characterization to prove $(c)\Rightarrow(a)$. Let $\alpha:I\rightarrow M$ be a curve starting at $p\in M$. By orthonormal expansion it suffices to prove that $R(X,Y,Z,W)$ is constant along $\alpha$, whenever $X,Y,Z,W$ are parallel vector fields along $\alpha$. Fix $t\in I$ and define a function $A: (T_pM)^4\rightarrow R$ by $$A(x,y,z,w)=R(X,Y,Z,W)(t),$$ where $X,Y,Z,W$ are parallel vector fields along $\alpha$ extending $x,y,z,w$, respectively. Given any $u\in T_pM$, if $U$ denotes its extension along $\gamma$, then $JU$ is the extension of $Ju$. Since $g((X\wedge_g Y)Y,X)$ is constant along $\alpha$, it holds $$\dfrac{A(u,Ju,Ju,u)}{g((u\wedge_g Ju)Ju,u)}=K(U,JU,JU,U)(t)=K(u,Ju,Ju,u).$$ It is easy to show that $A$ is curvaturelike (i.e., it is a function with the same properties as $R(x,y,z,w)$). Therefore, applying~\cite[Prop. 7.1, Ch. IX]{KN_1969}, we get $A=R$. Thus, $R(X,Y,Z,W)(t)$ is independent of $t$.
\end{demo}

The integrability condition of $\nabla R=0$ is $R\cdot R=0$. In the same manner as for locally symmetric Kaehler spaces, Kaehler semi-symmetric spaces can be characterized in terms of holomorphic planes, as it is shown in the next characterization result.

\begin{thm}\label{thm:semi-charac}
	Let $(M,g,J)$ be a Kaehler manifold. Then, the following conditions are equivalent:
	\begin{itemize}
		\item[$(a)$] $M$ is semisymmetric.
		\item[$(b)$] $R\cdot R(u,v,v,u;x,Jx)=0$ for every $p\in M$ and $x,u,v\in T_pM$.
		\item[$(c)$] $R\cdot R(u,Ju,Ju,u;x,y)=0$ for every $p\in M$ and $x,y,u\in T_pM$.
		\item[$(d)$] $R\cdot R(u,Ju,Ju,u;x,Jx)=0$ for every $p\in M$ and $x,u\in T_pM$.
	\end{itemize}
\end{thm}

\begin{demo}
	The proof is immediate. In fact, implications $(a) \Rightarrow (b) \Rightarrow (d)$ and $(a) \Rightarrow (c) \Rightarrow (d)$ are direct. It only remains to proof $(d) \Rightarrow (a)$, but it follows immediately from Proposition~\ref{prop:auxalg}.
\end{demo}

\section{Holomorphically pseudosymmetric Kaehler manifolds}\label{sec:holom_ps}

It is known that every Kaehler manifold $(M^{2n},g,J)$ of real dimension $2n>4$ which is pseudosymmetric in the sense of Deszcz is also semisymmetric, see~\cite{DDV} and~\cite{Olszak2003}. However, there exists an example in dimension $2n=4$ of a Kaehler manifold $(M^4,g,J)$ which is pseudosymmetric in the sense of Deszcz but not semisymmetric,~\cite{Olszak2003}. On the other hand, as it has been remarked in the introduction, the tensors $R\cdot R$ and $Q(g,R)$ do not present the same symmetries and properties involving the complex structure $J$. For these reasons, for Kaehler manifolds it seems reasonable to consider the alternative version of pseudosymmetry based on the complex Tachibana tensor. This definition was proposed by Olszak in 1989~\cite{Olszak89}.

\begin{definition}\label{def:pshol}
	A Kaehler manifold $(M^{2n},g,J)$, $n\geq 2$, is said to be holomorphically pseudosymmetric when the tensor $R\cdot R$ satisfies 
	\begin{equation}
	\label{eq:ps}
	R\cdot R=f\,\, Q^c(g,R),
	\end{equation}
	where $f\in\mathcal{C}^\infty(M)$. 
\end{definition}

It is obvious that any semisymmetric Kaehler manifold is always  holomorphically pseudosymmetric. However, the reverse is not true,  since examples of both compact and non-compact holomorphically pseudosymmetric Kaehler manifolds which are not semisymmetric are known, see~\cite{Jelonek} and~\cite{Olszak2009}, respectively.

Although the tensors involved in both definitions of pseudosymmetry are different, it is possible to characterize the holomorphic pseudosymmetry in terms of the double sectional curvatures as defined in~\eqref{eq:dousec}, getting results similar to those known for Riemannian manifolds,~\cite[Theorems 3 and 5]{HV_2007}. Using Proposition \ref{prop:auxalg}, we are now able to present our first result in this section.

\begin{thm}\label{thm:DvsO}
	A Kaehler manifold $(M^{2n},g,J)$, $n\geq 2$, is holomorphically pseudosymmetric if and only if for every $p\in\mathcal{U}$ and every holomorphic planes  in $T_pM$,  $\pi^h$ and $\bar{\pi}^h$, such that $\pi^h$ is curvature-dependent with respect to $\bar{\pi}^h$, the double sectional curvature $L(p,\pi^h,\bar{\pi}^h)$ is independent of both planes, i.e. $L(p,\pi^h,\bar{\pi}^h)=L(p)\in\mathcal{C}^\infty(U)$.
\end{thm}

\begin{demo}
	Let us assume first that $M$ is holomorphically pseudosymmetric, so in particular
	\begin{equation}\label{eq:holps_hplanes}
	R\cdot R(X,JX,JX,X;U,JU)=f\,\,Q^c(g,R)(X,JX,JX,X;U,JU)
	\end{equation}
	for every $X,U\in\mathfrak{X}(M)$, $f\in\mathcal{C}^\infty (M)$ being a smooth function. 
	
	Therefore, taking into account~\eqref{eq:TQ1} and the definition of the double sectional curvature, given $p\in \mathcal{U}$ and any pair of holomorphic planes in $T_pM$, $\pi^h=x\wedge Jx$ and $\bar{\pi}=u\wedge Ju$, such that $\pi^h$ is curvature-dependent with respect to $\bar{\pi}^h$, it holds $L(p,\pi^h,\bar{\pi}^h)=2\,f(p)$, so the double sectional curvature is independent of the planes.
	
	Conversely, let us suppose now that for any $p\in \mathcal{U}$ and any pair of holomorphic planes in $T_pM$, $\pi^h=x\wedge Jx$ and $\bar{\pi}^h=u\wedge Ju$, such that $\pi^h$ is curvature-dependent with respect to $\bar{\pi}^h$, its double sectional curvature does not depend of such planes, i.e.
	\begin{equation}\label{eq:lp_hplanes}
	R\cdot R(x,Jx,Jx,x;u,Ju)=L(p)\,\,Q(g,R)(x,Jx,Jx,x;u,Ju),
	\end{equation} 
	which taking into account~\eqref{eq:TQ1} implies
	\begin{equation}\label{eq:lp_hplanes_2}
	R\cdot R(x,Jx,Jx,x;u,Ju)=f(p)\,\,Q^c(g,R) R(x,Jx,Jx,x;u,Ju),
	\end{equation}
	where $f=\frac{1}{2}\,L\in\mathcal{C}^\infty(\mathcal{U})$.
	
	We claim that the same equality holds for any pair of, non necessarily curvature-dependent, holomorphic planes in $T_pM$, for all $p\in M$, for any smooth extension of $f$ in $M$. Therefore, since both $R\cdot R$ and $Q^c(g,R)$ satisfy the simmetries $(a)$, $(b)$, $(c)$ and $(d)$ of Proposition~\ref{prop:auxalg}, it follows that $R\cdot R=f\,\,Q^c(g,R)$, so $M$ is holomorphically pseudosymmetric.
	
	It remains to prove the claim. On the one hand, in the case $p\in\mathcal{U}$, let us consider the smooth function $q\in\mathcal{C}^\infty(T_pM\times T_pM)$, defined by $q(v,w)=Q(g,R)(v,Jv,Jv,v;w,Jw)$. It is easy to prove that the zero set of $q(v,w)$ does not contain any open subset. Fix a basis $\{e_i, Je_i:1\leq i\leq m\}$ and realize that $q(v,w)$ is a non null polynomial on the components of $v$ and $w$. If $\pi^h=x\wedge Jx\in T_pM$ is not curvature-dependent with respect to $\bar{\pi}^h=u\wedge Ju\in T_pM$, we have that $q(x,u)=0$. It is possible to choose sequences of tangent vectors in $T_pM$, $\left\{x_n\right\}_n$ and $\left\{u_n\right\}_n$, convergent to $x$ and $u$ respectively, such that $q(x_n,u_n)\neq 0$, $\forall n\in\mathbb{N}$. Equivalently, $x_n\wedge Jx_n$ is a curvature-dependent plane with respect to $u_n\wedge Ju_n$, for every $n\in\mathbb{N}$. Consequently, it holds
	\begin{equation}\label{eq:lp_hplanes_n}
	R\cdot R(x_n,Jx_n,Jx_n,x_n;u_n,Ju_n)=\frac{1}{2}\,L(p)\,Q^c(g,R)(x_n,Jx_n,Jx_n,x_n;u_n,Ju_n),
	\end{equation} 
	for all $n\in\mathbb{N}$, so~\eqref{eq:lp_hplanes_2} follows at $p\in\mathcal{U}$ by a continuity argument.
	
	On the other hand, the Tachibana tensor vanishes identically on $\textrm{int}(M\setminus \mathcal{U})$, so this open subset has constant sectional curvature and consequently~\eqref{eq:lp_hplanes_2} is trivially satisfied for any $f¨$. Finally, the equality also holds on $\partial(M\setminus\mathcal{U})$ by continuity, and the claim is proved.
	
\end{demo}

As an inmmediate consequence of Theorem~\ref{thm:DvsO} we get the following result.
\begin{corollary}\label{cor:implication}
	Any Kaehler manifold $(M^{2n},g,J)$ which is pseudosymmetric in the sense of Deszcz is also holomorphically pseudosymmetric.
\end{corollary}

The already mentioned examples in~\cite{Jelonek,Olszak2009} show that the converse is not true in general.

Finally, we can prove that it is possible to weaken the condition on the double sectional curvatures in the characterization of the holomorphic pseudosymmetry given in Theorem~\ref{thm:DvsO}. Specifically, to guarantee the holomorphic pseudosymmetry, it is only necessary to ask the double sectional curvatures to be independent of the first holomorphic plane.

\begin{thm}\label{thm:DvsO_weak}
A Kaehler manifold $(M^{2n},g,J)$, $n\geq 2$, is holomorphically pseudosymmetric if and only if for every $p\in\mathcal{U}$ and every holomorphic planes  in $T_pM$,  $\pi^h$ and $\bar{\pi}^h$, such that $\pi^h$ is curvature-dependent with respect to $\bar{\pi}^h$, the double sectional curvature $L(p,\pi^h,\bar{\pi}^h)$ is independent of $\pi^h$, i.e. $L(p,\pi^h,\bar{\pi}^h)=L(p,\bar{\pi}^h)$. 
\end{thm}

\begin{demo}
	According to Theorem~\ref{thm:DvsO} we only need to prove that if  all the double sectional curvatures at $p$, $L(p,\pi^h,\bar{\pi}^h)$, are independent of $\pi^h$, then they are independent of both planes, i.e.  $L\in\mathcal{C}^\infty(\mathcal{U})$.
	
	Following an analogous reasoning as the one in Theorem~\ref{thm:DvsO}, given $p\in \mathcal{U}$ and any pair of holomorphic planes in $T_p M$, $\pi^h=x\wedge Jx$ and $\bar{\pi}^h=u\wedge Ju$, such that $\pi^h$ is curvature-dependent with respect to $\bar{\pi}^h$, if $L(p,\pi^h,\bar{\pi}^h)=L(p,\bar{\pi}^h)$, then 
	\begin{equation}\label{eq:Dvs0_weak_1}
		R\cdot R(x_1,x_2,x_3,x_4;u,Ju)=\frac{1}{2}\,L(p,\bar{\pi}^h)\,Q^c(g,R)(x_1,x_2,x_3,x_4;u,Ju),
	\end{equation}
for every $x_1,x_2,x_3,x_4\in T_p M$. In particular,
\begin{equation}\label{eq:Dvs0_weak_2}
	R\cdot R(x,Jx,Jx,x;u,Ju)=\frac{1}{2}\,L(p,\bar{\pi}^h)\,Q^c(g,R)(x,Jx,Jx,x;u,Ju).
\end{equation}

Taking into account the algebraic symmetries of $R\cdot R$, see~\cite[Lemma 1]{HV_2007}, we get
\begin{equation}\label{eq:Dvs0_weak_symRR}
	\begin{split}
		R\cdot R(x,Jx,Jx,x;u,Ju)&=-R\cdot R(Jx,x,u,Ju;x,Jx)-R\cdot R(u,Ju,x,Jx;Jx,x)\\&=2\,R\cdot R(u,Ju,x,Jx;x,Jx).
	\end{split}
\end{equation}
Analogously, since $Q^c(g,R)$ satisfies the same algebraic properties as $R\cdot R$ it also yields
\begin{equation}\label{eq:Dvs0_weak_symPiR}
	Q^c(g,R)(x,Jx,Jx,x;u,Ju)=2\,Q^c(g,R)(u,Ju,x,Jx;x,Jx)\neq 0.
\end{equation}
Thus, from~\eqref{eq:Dvs0_weak_symRR} and~\eqref{eq:Dvs0_weak_symPiR}, expression~\eqref{eq:Dvs0_weak_2} becomes
\begin{equation}\label{eq:Dvs0_weak_3}
	R\cdot R(u,Ju,x,Jx;x,Jx)=\frac{1}{2}\,L(p,\bar{\pi}^h)\,Q^c(g,R)(u,Ju,x,Jx;x,Jx).
\end{equation}

We may observe at this point that the fact that the plane $\pi^h=x\wedge Jx$ is curvature-dependent with respect to $\bar{\pi}^h=u\wedge Ju$ implies that there exists necessarily an holomorphic plane, $\tilde{\pi}^h$, which is curvature-dependent with respect to $\pi^h$. Otherwise, Proposition~\ref{prop:auxalg2} would yield that $Q^c(g,R)(\cdot,\cdot,\cdot,\cdot;x,Jx)=0$, getting a contradiction of~\eqref{eq:Dvs0_weak_symPiR}. Therefore, we can reproduce the reasoning above applied to $\tilde{\pi}^h$ and $\pi^h$, obtaining \begin{equation}
	R\cdot R(u,Ju,x,Jx;x,Jx)=\frac{1}{2}\,L(p,\pi^h)\,Q^c(g,R)(u,Ju,x,Jx;x,Jx).
\end{equation}

From~\eqref{eq:TQ1}, neither $Q^c(g,R)(\pi^h;\bar{\pi}^h)$ nor $Q^c(g,R)(\tilde{\pi}^h;\pi^h)$ vanishes. Consequently, we conclude that $L(p,\pi^h)=L(p,\bar{\pi}^h)$. 

It remains to show that such equality is also satisfied for any pair of planes in $T_pM$, $p\in\mathcal{U}$, even if neither of them is curvature-dependent with respect to the other. To that end, let us asssume now that $\pi^h=x\wedge Jx$ is not a curvature-dependent plane with respect to $\bar{\pi}^h=u \wedge Ju$, nor vice versa. We can assume that there exists at least two holomorphic planes, $\pi_1^h=v_1\wedge Jv_1$ and $\pi_2^h=v_2\wedge Jv_2$, which are curvature-dependent with respect to $\pi^h$ and $\bar{\pi}^h$, respectively. Otherwise, there would be nothing to prove. Two options are possible. On the one hand, if $\pi_1^h$ is curvature-dependent with respect to $\pi_2^h$, or vice versa, then
\begin{equation}\label{eq:chainL}
L(p,\pi^h)=L(p,\pi_1^h)=L(p,\pi_2^h)=L(p,\bar{\pi}^h).
\end{equation}

On the other hand, if $\pi_1^h$ and $\pi_2^h$ are not curvature-dependent, it is always possible to construct another holomorphic plane $\pi_3^h=v_3\wedge Jv_3$ curvature-dependent with respect to both $\pi^h$ and $\bar{\pi}^h$, so the conclusion holds with a similar chain of equalities. In order to construct $\pi_3^h$, we just have to consider $v_3=\alpha_1v_1+\alpha_2v_2$ for certain real constants $\alpha_1$ and $\alpha_2$. Since $Q(g,R)(v_1,Jv_1,Jv_1,v_1;x,Jx)\neq 0$ and $Q(g,R)(v_2,Jv_2,Jv_2,v_2;u,Ju)\neq 0$, both functions $Q(g,R)(v_3,Jv_3,Jv_3,v_3;x,Jx)$ and $Q(g,R)(v_3,Jv_3,Jv_3,v_3;u,Ju)$ are non zero polynomials on $\alpha_1$ and $\alpha_2$, so there should be at least a choice for $\alpha_1$ and $\alpha_2$ for which both polynomials do not vanish simultaneously. 
\end{demo}

\section*{Acknowledgements}
The first and the third authors are partially supported by Andalusian FEDER 1380930-F.

\bibliographystyle{amsplain}

\end{document}